\documentclass[reqno,10pt, centertags,draft]{article}

\usepackage{amssymb}
\usepackage{amsthm}
\usepackage{amsmath}
\usepackage[mathscr]{eucal}
\usepackage{latexsym}

\numberwithin{equation}{section}

\newtheorem{theorem}{Theorem}[section]
\newtheorem{thm}{Theorem}[section]
\newtheorem{lemma}[thm]{Lemma}

\theoremstyle{definition}

\newtheorem{remark}[thm]{Remark}

\renewcommand{\Re}{{\rm Re}\,}
\renewcommand{\Im}{{\rm Im}\,}
\newcommand{\C}{\mathbb{C}}
\newcommand{\R}{\mathbb{R}}
\newcommand{\Z}{\mathbb{Z}}
\newcommand{\N}{\mathbb{N}}

\date{\today}

\begin{document}
\title {Eigenfunction expansions associated with 1d periodic 
differential operators of order $2n$}
\author{Vadim Tkachenko}
\maketitle

\begin{abstract}

\vspace{5pt}
We prove an explicit formula for the spectral expansions in $L^2(\R)$ generated by
selfadjoint differential operators 
$$
(-1)^n\frac{d^{2n}}{dx^{2n}}+\sum\limits_{j=0}^{n-1}\frac{d^{j}}{dx^{j}}
p_j(x)\frac{d^{j}}{dx^{j}},\quad p_j(x+\pi)=p_j(x),\quad x\in\R.
$$
\vspace{10pt}

\end{abstract}

\section{Statement of results} 

It is well known \cite{TIT1}, see also \cite{TIT2},
that for every Hill operator
\begin{equation}
H=- \frac{d^2}{dx^2} + q(x), \quad q(x)=q(x+\pi), \quad x\in\mathbb{R} \,
\label{1.1}
\end{equation}
with a real-valued potential function $q(x)$ there exists a sequence of real
numbers
$$
\mu_0=\mu_0^-=\mu_0^+< \mu_1^-\leq\mu_1^+<...<\mu_k^-
\leq\mu_k^+<...
$$
such that the spectrum of $L$ in the space $L^2(\mathbb{R})$
has the form 
$$
\sigma(L)=\bigcup\limits_{k=0}^\infty\;[\mu_k^+,\mu_{k+1}^-].
$$
In 1950 Titchmarsch proved (\cite{TIT1}, see also \cite{TIT2}) that 
every real-valued function $f\in L^2(\mathbb{R})$ may be represented 
in the form
\begin{align}\begin{split}
f(x)&=\frac{1}{\pi}\sum\limits_{n=0}^\infty\;(-1)^n\int\limits_{\mu_n^+}
^{\mu_{n+1}^-}{d\mu}\;p(\mu)
\left\{\phi(\pi,\mu)\theta(x,\mu)g(\mu) 
-\theta'(\pi,\mu)\phi(x,\mu)h(\mu)\right.\\
 &+\frac{1}{2}(\phi'(\pi,\mu)-\theta(\pi,\mu))\theta(x,\mu)
h(\mu)+\frac{1}{2}(\phi'(\pi,\mu)-\theta(\pi,\mu))\phi(x,\mu)
g(\mu)\},\label{1.2}
\end{split}\end{align}
where $\theta(x,\mu)$ and $\phi(x,\mu)$ are solutions of the equation 
$(H-\mu I)=0$ satisfying the initial conditions 
$\theta(0,\mu)=\phi'(0,\mu)=1,\;\theta'(0,\mu)=\phi(0,\mu)=0$, with
$p(\mu)=(4-(\theta(\pi,\mu)+\phi'(\pi,\mu))^2)^{-1/2}$ and
$$
g(\mu)=\int\limits_{\mathbb R}\;dx\;\theta(x,\mu)f(x),\quad
h(\mu)=\int\limits_{\mathbb R}\;dx\;\phi(x,\mu)f(x).
$$
The expansion \eqref{1.2} determines explicitly the spectral matrix of 
operator $H$ and, in particular, shows that the multiplicity of its spectrum 
equals $2$.

We consider arbitrary periodic self-adjoint differential operators
\begin{equation}
L=(-1)^n\frac{d^{2n}}{dx^{2n}}+\sum\limits_{j=0}^{n-1}\frac{d^{j}}{dx^{j}}
p_j(x)\frac{d^{j}}{dx^{j}},\quad p_j(x+\pi)=p_j(x),\quad x\in\R,\label{1.3}
\end{equation}
with real-valued functions $p_j(x),$ \ $j=0,1,...,n-1,$ such that 
\begin{equation}\label{1.31}
P(L)=\sum\limits_{j=0}^{n-1}\int_0^\pi|p^{(j)}_j(x)|\;dx<\infty.
\end{equation}

Similar to Hill operators, the spectrum of every such operator 
in the complex space $L^2(\R)$ has the band structure, 
but in contrast to these operators its multiplicity may vary inside one spectral 
band. Using a general resolvent method due to Kodaira and Spencer, Dunford and Scwartz
(cf., \cite{DAS}, Chap. XIII) proved a formula for the spectral matrix of operator 
\eqref{1.3} on an interval of a constant multiplicity of the spectrum.

Our aim is to obtain an expansion formula 
similar to \eqref{1.2} for operator \eqref{1.3}, to derive from it an explicit formula for the spectral matrix of 
operator \eqref{1.3} and to prove that this matrix determines its coefficients uniquely. 
Our approach is based 
on a version of the Fourier transform proposed by Gel'fand \cite{GEL}
for a study of periodic differential operators.

To state our main result, let $\{u_k(x,\mu)\}_{k=1}^{2n}$ be the 
fundamental system of solutions of equation
\begin{equation}
Ly-\mu y=0\label{1.4}
\end{equation}
normalized by the initial conditions
\begin{equation}
u_k^{(j)}(0,\mu)=\delta_{k-1,j},\quad k=1,...,2n;\quad j=0,...,2n-1,
\label{1.5}
\end{equation}
and let 
\begin{equation}\label{1.15}
U(\mu)=||u_k^{(j-1)}(\pi,\mu)||^{2n}_{k,j=1}
\end{equation}
be the monodromy matrix of $L.$ The eigenvalues of $U(\mu)$ are solutions of the 
characteristic equation 
\begin{equation}
\Delta(\mu,\rho)=0\label{1.6}
\end{equation}where $\Delta(\mu,\rho)=\det(U(\mu)-\rho I)$
and are called the Floquet multipliers of $L$. If $\{v_1,...,v_{2n}\}^\perp$ is an 
eigenvector corresponding to $\rho$, then the solution of \eqref{1.4} uniquely determined 
by the initial conditions
$$
y^{(j)}(0)=v_{j+1},\quad j=0,...,2n-1,
$$
has a ``quasi-periodic'' property
$$
y(x+k\pi)=\rho^ky(x),\quad x\in\R,\quad k\in\Z,
$$
and is called a Floquet solution of \eqref{1.4}.

The following proposition is well known in the theory of ordinary differential operators
with periodic coefficients, cf. \cite{DAS}.

\begin{theorem}\label{t1}
The spectrum $\sigma(L)$ of an operator \eqref{1.3} is absolutely continuous and 
coincides with the set of all 
those $\mu$ for which there exists a solution $\rho$ of \eqref{1.6} with
$|\rho|=1$.
\end{theorem}

Let us define the functions
\begin{equation}\label{1.12}
E(x;\mu,\rho)=\left|
\begin{array}{ccccc}
u_1(x,\mu)&...&u_j(x,\mu)&...&u_{2n}(x,\mu)\\
u_1(\pi,\mu)-\rho&...&u_j(\pi,\mu)&...&u_{2n}(\pi,\mu)\\
...&...&...&...&...\\
u_1^{(j-1)}(\pi,\mu)&...&u_j^{(j-1)}(\pi,\mu)-\rho&...&u_{2n}^{(j-1)}(\pi,\mu)\\
...&...&...&...&...\\
u_1^{(2n-2)}(\pi,\mu)&...&u_j^{(2n-2)}(\pi,\mu)&...&u_{2n}^{(2n-2)}(\pi,\mu)
\end{array}\right|
\end{equation}
and
\begin{equation}\label{1.11}
p(\mu,\rho)=|2\pi\; E(0;\mu,\rho^{-1})\Delta'_\rho(\mu,\rho)|^{-1}.
\end{equation}

We will see later on that with a proper choice of $t$ and $\mu$ the functions  
$E(x;\mu,e^{it})$ are the Floquet solutions participating in the spectral 
expansion associated 
with $L$. As to the function \eqref{1.11}, it supplies the normalizing factors 
in such expansion.

\begin{theorem}\label{t6} If $\omega_k,\: k=1,...,2n,$ are all values of 
$\sqrt[2n]{(-1)^n}$ then  for every $k=1,...,2n$
there exists a solution $\rho_k(\mu)$ of \eqref{1.6}, continuous on the 
real line and satisfying the asymptotic relation
\begin{equation}
|\rho_k(\mu)|=e^{\Re\omega_k\lambda\pi}(1+o(1)),\quad\mu=\lambda^{2n},
\quad\lambda\to+\infty.\label{2.20}
\end{equation}
Moreover, solutions $\rho_k(\mu)$ are pair-wise distinct and 
analytic at points of $\R$, except maybe points of a discrete 
set grouped in pairs asymptotically close to the set 
\begin{equation}\label{1.16}
{\cal N}=\bigcup\limits_{j=1}^{2n}N_j,\quad {\cal N}_j=
\{(-1)^{j}m^{2n}(\Im\omega_j)^{-2n}\}_{m=1}^\infty
\end{equation}
where some of them coincide and their analyticity may fail.
\end{theorem}
If we set 
\begin{equation}\label{1.13}
\sigma_k(L)=\{\mu\in\R:|\rho_k(\mu)|=1\},\qquad k=1,...2n,
\end{equation}
then we obtain 
$$
\sigma(L)=\bigcup\limits_{k=1}^{2n}\;\sigma_k(L).
$$

Denote by ${\mathcal H}^2_{2n}(L)$ the Hilbert space of complex 
$2n$-vector 
functions $\Phi=
\{\phi_k(\mu)\}_{k=1}^{2n}$ on the sets $\{\sigma_k(L)\}_{k=1}^{2n}$ with the 
scalar product
$$
(\Phi,\Psi)=\sum\limits_{k=1}^{2n}
\int\limits_{\sigma(L)}\;d\mu\;\chi_k(\mu)\;p(\mu,\rho_k(\mu))\;
\phi_k(\mu)\;\overline{\psi_k(\mu)}
$$
where $\chi_k(\mu)$ is the indicator function of the set  ${\sigma}_k(L)$.

For the simplest operator 
\begin{equation}\label{1.14}
(-1)^n\;\frac{d^{2n}}{dx^{2n}}
\end{equation} 
we have $\rho_k(\mu)=\exp(\omega_k\mu^{1/2n}\pi)$
and the sets $\sigma_k(L)$ are reduced to the point $0$ for all $k$'s 
with $\omega_k\neq\pm i.$ Therefore a 
situation where some intervals $\sigma_k(L)$ degenerate to a point cannot be ruled out.
As a result, some functions 
$\phi_k(\mu)$ may be trivial for all $\Phi\in{\mathcal H}^2_{2n}(L)$.
\begin{theorem}\label{t2}
The relations 
\begin{equation}\label{1.7}
\Phi(\mu,\rho_k(\mu);f)=\int\limits_\R\;dy\;f(y)\;E(y;\mu,(\rho_k(\mu))^{-1}),
\qquad k=1,...,2n,
\end{equation}
and 
\begin{equation}\label{1.8}
f(x)=
\int\limits_{\sigma(L)}\;d\mu\;
\sum_{k=1}^{2n}\chi_k(\mu)\;p(\mu,\rho_k(\mu))\phi_k(\mu)\;E(x;\mu,\rho_k(\mu))
\end{equation}
define one-to-one mapping of $L^2(\R)$ onto ${\mathcal H}^2_{2n}(L)$ and its inverse 
conjugating 
operator $L$ on the former space with the scalar operator $\mu I$ on the latter space.

The integrals \eqref{1.7} and \eqref{1.8} converge in the norms of the corresponding spaces 
and for every function $f\in L^2(\R)$ the Parseval identity 
\begin{equation}\label{1.9}
\int\limits_\R\;dx\;|f(x)|^2=
\int\limits_{\sigma(L)}\;d\mu\;
\sum_{k=1}^{2n}\;\chi_k(\mu)p(\mu,\rho_k(\mu))\;|\Phi(\mu,\rho_k(\mu);f)|^2
\end{equation}
holds.
\end{theorem}

It is easy to see that in the case of Hill operators Theorem \ref{t2} states 
that for every function $f\in L^2(\R)$ the representation 
$$
f(x)=\frac{1}{4\pi}\int\limits_{\sigma(L)}\;d\mu\;\frac{|\phi(\pi,\mu)|}
{\sqrt{1-u_+(\mu)^2}}\left\{Y_+(x,\mu)F_-(\mu;f)+Y_-(x,\mu)F_+(\mu;f)
\right\}
$$
is valid where 
$$
Y_\pm(x,\mu)=\theta(x,\mu)-\frac{u_-(\mu)\pm i\sqrt{1-u_+(\mu)^2}}
{\phi(\pi,\mu)}\;\phi(x,\mu)
$$
are the Floquet solutions,
$$
F_\pm(\mu;f)=\int\limits_{\mathbb{R}}\;dy\;f(y)\;Y_\pm(y,\mu),
$$
and
$$
u_\pm(\mu)=\frac{\theta(\pi,\mu)\pm\phi'(\pi,\mu)}{2}.
$$
This is a complex version of expansion \eqref{1.2} obtained first in \cite{MEI} and 
later in \cite{PAT}.

\section{Solutions $\rho_k(\mu)$ and Riemann surface ${\cal R}(L)$}

Let us start from the proof of Theorem \ref{t6}. Following \cite{NAI} we enumarate 
the numbers $\omega_k=\sqrt[2n]{(-1)^n}$ in such an order that 
\begin{equation}\label{2.7}
\begin{array}{c}
1=\omega_{1}>\Re\omega_{2}=\Re\omega_{3}>...>\Re\omega_{2k}=\Re\omega_{2k+1}>...>\omega_
{2n}=-1,\\
\Im\omega_{2k}\geq0,
\end{array}
\end{equation}
if $n$ is even, and 
\begin{equation}\label{2.8}\begin{array}{c}
\Re\omega_{1}=\Re\omega_{2}>...>\Re\omega_{2k+1}=\Re\omega_{2k}>...>\Re\omega_
{2n-1}=\Re\omega_{2n},\\
\Im\omega_{2k}\leq0,
\end{array}
\end{equation}
if $n$ is odd. With such enumeration we have $\omega_n=i$ both for odd and even values 
of $n$. In what follows we restrict ourselves to the case $n=2p, p\in\N$. The alternative 
case differs from it by non-essential technical details.

Let 
$$T(r)=\{\lambda:\lambda=z-e^{-i\pi/4n},\quad|z|\geq r,\quad-\frac{1}{2n}\leq\arg z\leq0\}.$$
According to \cite{NAI} there exists a fundamental matrix $Y(x,\mu)$ of solutions
to \eqref{1.4} representable in the form
\begin{equation}\label{2.11}
Y(x,\mu)=D(\lambda)\Omega(x,\lambda)E(\lambda x),\quad\mu=\lambda^{2n},\quad x\in[0,\pi],
\quad\lambda\in T(r),
\end{equation}
with matrices
$$\begin{array}{cc}
D(\lambda)=\|\lambda^{j-1}\delta_{jk}\|_{j,k=1}^{2n},&
E(x)=\|e^{\omega_kx}\delta_{jk}\|_{j,k=1}^{2n},\\
 & \\
\Omega(x,\lambda)=\Omega+\lambda^{-1}\Omega_1(x,\lambda),&
\Omega=\|\omega_k^{j-1}\|_{j,k=1}^{2n},
\end{array}$$
where
$$\sup\limits_{
x\in[0,\pi],\;
\lambda\in T(r)}\|\Omega_1(x,\lambda)\|<\infty.
$$

Since $U(\pi,\mu)=Y(x,\mu)Y(0,\mu)^{-1}$, we can represent \eqref{1.6} in the form
\begin{equation}\label{2.4}
\det((I+o(1))E(\pi\lambda)-\rho I)=0.
\end{equation}

To investigate the latter equation we define the entire functions
$$
f_k(\lambda)=1-e^{(\omega_{k+1}-\omega_k)\lambda\pi},\qquad k=1,...,2n-1,
$$
and denote by $Z_k$ their zero sets. 
For every $\delta>0$ we denote by $U_k(\delta)$ the $\delta$-
neighborhood of $Z_k$. If  
$$
T_k(r,\delta)=\left\{
\begin{array}{cc}
T(r)\setminus U_1(\delta),&k=1,\\
 & \\T(r)\setminus (U_k(\delta)\bigcup U_{k-1}(\delta)),&1<k<2n-1,\\
 & \\T(r)\setminus U_{2n-1}(\delta),&k=2n,
\end{array}
\right.
$$ 
then 
\begin{equation}\label{2.21}
\begin{array}{cc}
C^{-1}\delta\leq |f_k(\lambda)|\leq C\delta,&\lambda\in\partial U_k(\delta)
\bigcap T_k(r,\delta),\\
 & \\C^{-1}\delta\leq |f_k(\lambda)|\leq C,&\lambda\in T_k(r,\delta).
\end{array}\end{equation}
with a constant $C>1$ not depending on $\delta$.
In what follows we fix $\delta$ such that $C\delta<10^{-1}$.
 
First, we the substitute $\rho=\sigma e^{\omega_1\lambda\pi}$
in \eqref{2.4}, divide the resulting equation by $\rho=e^{2n\omega_1\lambda\pi}$
and obtain the equivalent equation
\begin{equation}\label{2.5}
\sigma^{2n-2}(\sigma-1)(\sigma-1-f_2(\lambda))=Q_1(\sigma,\lambda),
\end{equation}
where $Q_1(\sigma,\lambda)$ is a polynomial in $\sigma$ 
with coefficients analytic in $T_1(r,\delta)$ and vanishing as $|\lambda|\to+\infty$. 
For sufficiently small $\epsilon>0$ we find large $r$ such that the function
$$
|Q_1(\sigma,\lambda)\sigma^{-(2n-2)}(\sigma-1-f_2(\lambda))^{-1}|\leq\frac{\epsilon}{2},
\quad\lambda\in 
T_1(r,\delta),\quad|\sigma-1|\leq\epsilon,$$ 
and using the Rouchet Theorem conclude that there exists the unique solution  
$\sigma_1(\lambda)$ of \eqref{2.5} analytic in $T_1(r,\delta)$ and such that 
$|\sigma_1(\lambda)-1|\leq\epsilon$. Since $(r,\infty)\in T_1(r,\delta)$, Equation
\eqref{2.5} implies that
$$
\rho_1(\mu)=e^{\omega_1\lambda\pi}\sigma_1(\lambda),\quad\mu=\lambda^{2n},\quad \lambda\in
T_1(r,\delta),
$$
is a solution of \eqref{1.6} satisfying \eqref{2.20} with $k=1$.

In the same way the substitution  $\rho=\sigma e^{\omega_{2n}\lambda\pi}$ leads us to the equation 
\begin{equation}\label{2.6}
(\sigma-1)(\sigma-1-\sigma f_{2n-1}(\lambda))=Q_{2n}(\sigma,\lambda),\quad \lambda\in
T_{2n}(r,\delta),
\end{equation}
with its unique solution $\sigma_{2n}(\lambda)$ analytic in $T_{2n}(r,\delta)$ and such that 
$$
\rho_{2n}(\mu)=e^{\omega_{2n}\lambda\pi}\sigma_{2n}(\lambda),\quad\mu=\lambda^{2n},\quad \lambda\in T_{2n}(r,\delta),
$$
is a solution of \eqref{2.4} satisfying \eqref{2.20} with $k=2n$.

If now $1<k=2p<2n$ and $\rho=\sigma e^{\omega_k\lambda\pi}$ then Equation \eqref{2.4} takes on the form
\begin{equation}\label{2.23}
\sigma^{2n-k-1}(\sigma-1-\sigma f_{k-1}(\lambda))(\sigma-1)(\sigma-1+f_{k}(\lambda))=
Q_{k}(\sigma,\lambda),\quad \lambda\in T_{k}(r,\delta),
\end{equation}
where $Q_{k}(\sigma,\lambda)$ is a polynomial in $\sigma$ 
with coefficients analytic in $T_k(r,\delta)$ and vanishing as $|\lambda|\to+\infty$.
As before we find that Equation \eqref{2.23} has the unique solution $\sigma_{k}(\lambda)$ analytic in $T_{k}(r,\delta)$ and satisfying the relation 
\begin{equation}\label{2.25}
\lim\limits_{\lambda\in T_{k}(r,\delta);|\lambda|\to\infty}\sigma_k(\lambda)=1
\end{equation}
but now, contrary to the case $k=1$ and $k=2n$, the set 
$T_{k}(r,\delta)$ does not contain a part of the ray $(r,+\infty)$ belonging to the set 
$U_k(\delta)$.

To define $\sigma_{k}(\lambda)$ inside the exceptional set $U_k(\delta)$ we assume, to 
be definite, that  
$\alpha\in{\cal N}_k\bigcap T(r)$, define $D_\alpha=\{\lambda:|\lambda-\alpha|\leq\delta\}$ and represent \eqref{2.23} in the form
\begin{equation}\label{2.24}
(\sigma-1)(\sigma-1+f_{k}(\lambda))=q_k(\sigma,\lambda),\quad \lambda\in D_\alpha,\quad
\end{equation}
where
$$
q_k(\sigma,\lambda)=Q_{k}(\sigma,\lambda)\sigma^{-2n+k+1}(\sigma-1-\sigma f_{k-1}(\lambda))^{-1}.
$$
If necessary we increase $r$ to satisfy
the estimate
$$
|q_{k}(\sigma,\lambda)|\leq 2C\delta^2,\quad|\sigma-1|\leq C^2,\quad\lambda\in D_\alpha,
$$
where $C$ is the constant from \eqref{2.21}. According to the Rouchet Theorem for every $\lambda\in D_\alpha$ there exist two solutions of \eqref{2.24} satisfying $|\sigma-1|\leq2C\delta$. If these solutions coincide at some point $\lambda\in D_\alpha$ then 
$f_{k}^2(\lambda)=4q_k(\sigma,\lambda)$. Once again we use the Rouchet Theorem and find that 
there are at most two such points $\lambda$. The analyticity of $q_k(\sigma,\lambda)$ in the domain $\{(\sigma,\lambda):|\sigma-1|\leq2C\delta, \lambda\in D_\alpha\}$ implies that there exists two-valued analytic solution $\tilde{\sigma}_k(\lambda)$ of \eqref{2.24} in $D_\alpha$ satisfying the estimate 
\begin{equation}\label{2.27}
|\tilde{\sigma}_k(\lambda)-1|\leq2C\delta,\qquad \lambda\in D_\alpha.
\end{equation} 
The function ${\sigma}_k(\lambda),\;\lambda\in T_k(r,\delta),$ is a single-valued solution
of Eq.\eqref{2.23} which is equivalent to Eq.\eqref{2.24} in $D_\alpha$ and because of 
\eqref{2.25} it coincides with a branch of $\tilde{\sigma}_k(\lambda)$. In other words, 
$\sigma_k (\lambda)$ is extended as a two-valued analytic function inside $D_\alpha$ with 
at most two ramification points.

Furthermore, with the same value of even $k=2p$ we substitute 
$\rho=\sigma e^{\omega_{k+1}\lambda\pi}$ in \eqref{2.4} and obtain an equation of the 
same type as \eqref{2.23}.
  As before we prove that there exists its unique solution 
$\sigma_{k+1}(\lambda)$
analytic in $T_{k+1}(r,\delta)$ and satisfying 
\begin{equation}\label{2.28}
\lim\limits_{\lambda\in T_{k+1}(r,\delta);|\lambda|\to\infty}\sigma_{k+1}(\lambda)=1.
\end{equation}
Similar to $\tilde{\sigma}_k(\lambda)$ its two-valued analytic extension $\tilde{\sigma}_{k+1}(\lambda)$ inside exceptional discs $D_\alpha$ satisfies condition 
\begin{equation}\label{2.29}
|\tilde{\sigma}_{k+1}(\lambda)-1|\leq2C\delta,\qquad \lambda\in D_\alpha.
 \end{equation} 

Let us now consider two-valued solutions  
$$
\tilde{\rho}_k(\mu)=\tilde{\sigma}_k(\lambda) e^{\omega_{k}\lambda\pi},\quad
\tilde{\rho}_{k+1}(\mu)=\tilde{\sigma}_{k+1}(\lambda) e^{\omega_{k+1}\lambda\pi},\quad
\mu=\lambda^{2n}, \quad\lambda\in D_\alpha,
$$
of \eqref{1.8}. The coefficients of Equation \eqref{1.6} are real on the real line and 
therefore the function
$$
\rho_k^*(\mu)=\overline{\rho_k(\overline{\mu})},\qquad \mu=\lambda^{2n},\qquad\lambda
\in\partial D_\alpha,
$$
is its solution as well. Since $\overline{\omega}_k={\omega}_{k+1}$, we have 
\begin{equation}\label{2.30}
\rho_k^*(\mu)={\sigma}_k^*({\lambda})e^{\omega_{k+1}\lambda\pi}=\rho_{k+1}(\mu),
\qquad \mu=\lambda^{2n},\qquad\lambda
\in\partial D_\alpha.
\end{equation}
On the other hand, 
$$\rho_k^*(\mu)={\theta}_k({\lambda})e^{\omega_{k}\lambda\pi},\quad
{\theta}_k({\lambda})={\sigma}_k^*(\lambda)e^{(\omega_{k+1}-\omega_{k})\lambda\pi},
\qquad \mu=\lambda^{2n},\quad\lambda
\in\partial D_\alpha,
$$
and the estimate 
$$
|{\theta}_k({\lambda})-1|\leq |\sigma_k(\lambda)f_k(\lambda)|+|\sigma_k(\lambda)-1|
\leq2C\delta,\qquad\lambda\in\partial D_\alpha,
$$
shows that $\rho_k^*(\mu)$ is a single-valued analytic solution of \eqref{2.24} in $\partial D_\alpha$.
If $\rho_k(\mu)=\rho_k^*(\mu)$ for $\mu=\lambda^{2n}, \lambda\in\partial D_\alpha$ then 
$|f_k(\lambda)|=o(1)$ which is impossible for large $\lambda$. Therefore for $\mu=\lambda^{2n}, \lambda\in\partial D_\alpha$, the functions $\rho_k(\mu)$ and 
$\rho_k^*(\mu)$ are different branches of the two-valued function $\tilde{\rho}_k(\mu)$.

To estimate the latter functions inside $D_\alpha$ we set
$$
A_k(\lambda)=\frac{\tilde{\rho}_k(\lambda^{2n})+\tilde{\rho}_k^*(\lambda^{2n})}{2}e^{-\lambda\pi\Re\omega_k},\quad
\lambda\in D_\alpha,
$$
and
$$
B_k(\lambda)={\tilde{\rho}_k(\lambda^{2n})\tilde{\rho}_k^*(\lambda^{2n})}e^{-2\lambda\pi\Re\omega_k},\quad\lambda\in D_\alpha.
$$
According to their definition both functions are single-valued and analytic 
inside $D_\alpha$ and, as it follows from \eqref{2.25}, \eqref{2.28} and  \eqref{2.30}, 
satisfy conditions
$$A_k(\lambda)=
\cos(\lambda\pi \Im\omega_k)+o(1),\quad
B_k(\lambda)=
1+o(1), \quad\lambda\in\partial D_\alpha.
$$
According to the Maximum Principle the same representations are valid in the domain $D_\alpha$.
It means that  $\tilde{\rho}_k(\mu)e^{-\lambda\pi\Re\omega_k}$ and 
$\tilde{\rho}_{k+1}(\mu)e^{-\lambda\pi\Re\omega_{k+1}}$ 
are solutions of the equation $
w^2-2(\cos(\lambda\pi \Im\omega_k)+o(1))w+1+o(1)=0,\quad\lambda\in D_\alpha,$
and therefore $ \tilde{\rho}_k(\mu)=e^{\omega_k\lambda\pi}(1+o(1)),  
\tilde{\rho}_{k+1}(\mu)=e^{\omega_{k+1}\lambda\pi}(1+o(1)), \quad\lambda\in D_\alpha.$
To complete the proof of Theorem \ref{t6} we fix a point $\mu=\lambda^{2n}$ with
$\lambda\in T_{k}(r,\delta)$, a system $\rho_k(\mu)=e^{\omega_k\lambda\pi}\sigma_k(\lambda),
k=1,...,2n,$ of solutions of \eqref{1.6} and extend all of them to $\R$ as single-valued
continuous functions, pair-wise distinct outside the discriminant set
$$
{\cal Z}(L)=\{\mu:\Delta(\mu,\rho)=\Delta'_\rho(\mu,\rho)=0\}.
$$
The latter is a zero set of the resultant $R[\Delta,\Delta'_\rho]$, see \cite{VDW}.
Since there exist points $\mu$ at which Equation \eqref{1.6} has $2n$ distinct roots,
the resultant is a non-trivial function and the set ${\cal Z}(L)$ has no
finite accumulation points. As a result, the extended solutions are analytic outside 
${\cal Z}(L)$ and, according to the above estimates of $\tilde{\rho}_k(\mu)$ and 
$\tilde{\rho}_k^*(\mu)$, satisfy \eqref{2.20}, which completes the proof of Theorem \ref{t6}.


The entries of the monodromy matrix are entire functions of order $1/2n,$ and Equation
\eqref{1.6} defines the $2n$-sheeted Riemann surface
$$
{\mathcal R}(L)=\{(\mu,\rho):\Delta(\mu,\rho)=0\}
$$
with the analytic function $\rho(\mu)$ on it.
\begin{lemma}\label{l2.0}
The Riemann surface ${\mathcal R}(L)$ is simply-connected.
\end{lemma}

Proof. Let $\{\rho_k(\mu)\}_{k=1}^{2n}$ be the system 
of unique solutions of \eqref{1.6} obtained in the proof of Theorem \ref{t6}.
As we have seen, these solutions are single-valued and analytic in the domains
$$
\Pi_k^-(r,\delta)=\{\mu:\mu=\lambda^{2n},\;\lambda\in T_k(r,\delta)\}
\bigcap\{\mu:\Im\mu\leq0\},
$$ 
and satisfy the asymptotic relations
$$
\rho_k(\mu)=e^{\omega_k\lambda\pi}(1+o(1)),\quad\mu=\lambda^{2n},
\quad|\lambda|\to\infty,\quad\lambda\in T_k(r,\delta).
$$
We will use their analytic extensions to describe the surface ${\mathcal R}(L)$.

First, we set
$$
\Pi_k^+(r,\delta)=\{\mu: \mu=\overline{z},\;z\in\Pi_k^-(r,\delta)\}
$$
and note that the function 
$\rho_k^*(\mu)=\overline{\rho_k(\overline{\mu})}$
are single-valued and analytic solutions of \eqref{1.6} in the domain 
$\Pi_k^+(r,\delta)$.

Furthermore, let
$$
c_k^+=\partial\Pi_k^-(r,\delta)\bigcap\{\mu\in\R:\mu\geq0\},\;
c_k^-=\partial\Pi_k^-(r,\delta)\bigcap\{\mu\in\R:\mu\leq0\}.
$$
Since the number $\omega_1$ is real, the solutions $\rho_1(\mu)$ and $\rho_1^*(\mu)$
have the same asymptotic behavior in $c_1^+$ and hence they coincide in $c_1^+$. 
By the same 
reason the solutions $\rho_{2n}(\mu)$ and $\rho_{2n}^*(\mu)$ coincide in $c_{2n}^+$.

For $k=2p,\: p=1,...,n-1,$ the function $\rho_k(\mu)$ coincides with $\rho_{k+1}^*(\mu)$ 
on the set  $c_k^+=c_{k+1}^+$ and with $\rho_{k-1}^*(\mu)$ on the set  $c_k^-=c_{k-1}^-$,
while  $\rho_{k+1}(\mu)$ coincides with $\rho_{k}^*(\mu)$  
on the set  $c_{k+1}^+$ and with $\rho_{k+2}^*(\mu)$ on the set  $c_{k+1}^-=c_{k+2}^-$.
If we glue together the pairs of the corresponding sets belonging to $\Pi_k^-(r,\delta)$ 
and $\Pi_k^+(r,\delta)$, then we obtain a Riemann surface ${\cal R}(r,\delta)$ with 
the single-valued analytic function $\rho(\mu)$ on it. The surface ${\cal R}(r,\delta)$
 is the same 
for all operators of the form \eqref{1.3} and may be obtained  from the surface 
${\cal R}_0=\{(\mu, \lambda): \mu=\lambda^{2n}\}$ corresponding to the simplest operator 
\eqref{1.14} after removing from it the disc $\{\mu:|\mu|\leq r\}$ and small neighborhoods  
of points projecting into the set $Z$ from \eqref{1.16}. It is evident that the surface 
${\cal R}_0$ is simply-connected. The surface ${\cal R}(L)$ results 
from the analytic extension of all functions $\rho_k(\mu)$ inside all exceptional sets.

To prove Lemma \ref{l2.0} let us assume that $(\mu_0,\rho(\mu_0))$ is an
arbitrary non-ramified
 point of ${\cal R}(L)$ and $(\nu_0,\rho(\nu_0))$ is a point of the
surface ${\cal R}(r,\delta)$. Denote by $l$ a simple smooth curve in the complex
plain connecting $\mu_0$ with $\nu_0$ and not containing points of the
discriminant set ${\cal Z}(L)$. According to the Monodromy Theorem \cite{MRK}
there exists the unique analytic continuation of $\rho(\mu)$ from a neighborhood
of $\mu_0$ along $l$. The regular element of this continuation at a
neighborhood of $\nu_0$ is a locally single-valued analytic solution of  
\eqref{1.6} and since the system $\{\rho_k(\mu)\}
_{k=1}^{2n}$ contains all local solutions analytic at $\nu_0$, this element
 coincides with some function $\rho_k(\mu)$.
It means that the lifting of $l$ to the surface ${\cal R}(L)$ connects
the points $(\mu_0,\rho(\mu_0))$ and $(\nu_0,\rho_k(\nu_0))$ which proves
the lemma.

We can give now a geometric description of the surface ${\cal R}(L)$. 
To this aim denote  by ${\cal R}_k$ a copy of the complex plane cut along 
the following sets:

1. A simple smooth curve inside the disc $\{\mu:|\mu|\leq r\}$ containing the point $\mu=-r$ and all points of the discriminant set ${\cal Z}(L)$ lying 
inside it;

2. The ray $\{\mu:\mu\leq -r\}$;
 
3. All real segments $[\alpha,\beta]$ where  $\alpha$ and $\beta$ are neighboring
 ramification points of $\rho_k(\mu)$ with $\Re\mu\geq r$;

4. All segments  $[\alpha,\overline\alpha]$ where $\alpha$ is a ramification 
point of either $\rho_{k-1}^*(\mu)$ with  $\Re\mu\leq-r$ if $k$ is even or 
$\rho_{k+1}^*(\mu)$ with  $\Re\mu\leq-r$ if $k$ is odd.

The function $\rho_k(\mu)$ is extended from $\Pi_k^-(r,\delta)$ to  ${\cal R}_k$ 
as 
a single-valued analytic solution  of \eqref{1.6} and by glueing the sheets 
${\cal R}_j, j=1,...,2n,$ according to boundary values of functions 
 $\rho_j(\mu)$ we obtain the surface ${\cal R}(L)\}$ with a single valued analytic function $\rho(\mu)$ on it.

{\bf Corollary}. {\it The transformation $J(\mu,\rho(\mu))=(\mu,(\rho(\mu))^{-1})$ is an 
analytic involution in  ${\cal R}(L)$.}

Indeed, if we have $|\rho_n(\mu)|=1$ with a sufficiently large real $\mu$ then 
$(\rho_n(\mu))^{-1}=\overline{\rho_n(\mu)}=\rho_{n+1}(\mu)$. It means that 
$(\rho(\mu))^{-1}$ is a solution of \eqref{1.6}. 
Since $\rho(\mu)$ is an analytic function in a simply-connected Riemann 
surface ${\cal R}(L)$, the function $(\rho(\mu))^{-1}$ is extended as a solution 
to entire surface and $(\mu,(\rho(\mu))^{-1})$ is its point, which proves Corrolary.

To conclude the present section, we note that the characteristic polynomial 
of operator $L$ is of the form 
$$
\Delta(\mu,\rho)=\rho^{2n}+\sum_{k=1}^{2n-1}A_k(\mu)\rho^k+1
$$
where $A_k(\mu)$ are entire functions. It follows from Corollary that 
$A_k(\mu)=A_{2n-k}(\mu)$ for $k=1,...,2n-1$. A statement of such a type for 
canonical Hamiltonian systems is known as the Lyapunov-Poincar\'e Theorem,
cf., \cite{KRE,YAS}.

\section{Band structure of the spectrum of operator $L$}
To describe the band structure of the spectrum $\sigma(L)$ let us introduce the set
$$
{\mathcal E(L)}=\{\mu\in\R:\mbox{there exists}\; \rho\in\C,
|\rho|=1, \mbox{ such that}\; \Delta(\mu,\rho)=\Delta'_\rho(\mu,\rho)=0\}.
$$
If 
\begin{equation}\label{2.17}
{\mathcal E_k(L)}=\{\mu\in\sigma(L):\Delta(\mu,\rho_k(\mu))=
\Delta'_\rho(\mu,\rho_k(\mu))=0\}
\end{equation}
then
$$
{\mathcal E(L)}=\bigcup\limits_{k=1}^{2n}\;\mathcal{E}_k(L).
$$
The set ${\mathcal E(L)}$ is a part of the discriminant set ${\cal Z}(L)$ and hence it is countable with the unique accumulation point at $+\infty$.

\begin{lemma}\label{l2.2}
If  $\rho_k(\mu_0)=e^{it_0}$ for some real $\mu_0\notin{\mathcal E_k(L)}$, $t_0\in[0,\pi]$
and some integer $k, 1\leq k\leq 2n,$ then there exist the maximal closed interval 
$S\subset[0,\pi]$ containing $t_0$ and the continuous monotonic function $\mu(t)$ 
in ${S}$ such that

(\;i) The relations 
\begin{equation}\label{2.15}
\Delta(\mu(t),e^{it})=0,\quad t\in{S}; \quad\mu(t_0)=\mu_0;\quad\mu'(t)\neq0,\quad t\in
\mbox{int}\ S,
\end{equation}
are valid;

(ii) The function $\mu(t)$  
maps ${S}$ one-to-one onto the compact interval $\ell\subset\sigma_k(L)$ 
with end-points in the set ${\mathcal E_k(L)}$.
\end{lemma}

Proof. Since  $\mu_0\notin{\mathcal E_k(L)}$, the surface ${\mathcal R}(L)$ is not 
ramified at the point  $(\mu_0,\rho_k(\mu_0))$ and $\rho_k(\mu)$ is a single-valued 
branch of 
$\rho(\mu)$ analytic in a small complex neighborhood 
$V_0=\{\mu:|\mu-\mu_0|<\epsilon\}\subset{\mathcal R}_k$ of $\mu_0$. 
Its Taylor expansion at  $\mu_0$ has the form 
\begin{equation}\label{2.13}
\rho_k(\mu)=e^{it_0}+c_k(\mu-\mu_0)^p+\sum\limits_{m=p+1}^\infty\;c_{km}(\mu-\mu_0)^{m},
\quad c_k\neq0,\quad p\geq1.
\end{equation}
If here $p\geq2,$ then the pre-image with respect to $\rho_k(\mu)$ of a small 
neighborhood 
of the point $e^{it_0}$ on the unit circle $U_0=\{z:|z|=1\}$ contains non-real points.
 According to Theorem \ref{t1} such points belong to the spectrum $\sigma(L)$ 
which is impossible,
since $L$ is a selfadjoint operator. Therefore $p=1$ in \eqref{2.13}, $\rho_k'(\mu_0)=c_k
\neq0$ and the function $\rho_k(\mu)$ maps $V_0$ one-to-one onto a small complex 
neighborhood
$\Theta_0$ of the point  $e^{it_0}\in U_0.$\footnote{These arguments are similar to
those used in \cite{DAS}, Ch.XIII.}

According to Theorem \ref{t1} the pre-image of the arc $\Theta_0\cap U_0$ consists of
the points belonging to the spectrum $\sigma(L)$ and hence coincides with the interval 
$W_0=V_0\cap\R$, $\mu_0$ being its inner point. The function 
\begin{equation}\label{2.14}
\rho^*_k(\mu)=\overline{\rho_k(\overline{\mu})},\quad \mu\in V_0,
\end{equation}
is a single-valued branch of $\rho(\mu)$ in $V_0$ different from $\rho_k(\mu)$ and 
if  $\rho_k(\mu_0)=\pm1$ then  $\rho_k(\mu_0)=\rho^*_k(\mu_0)$ contradicting the 
assumption $\mu_0\notin{\mathcal E_k(L)} $. Therefore $t_0\notin\{0,\pi,2\pi\}$.

Denote by $\ell$ the largest closed interval in $\sigma_k(L)$ which contains $\mu_0,$ 
but does not contain points of the set ${\mathcal E}_k(L)$ in its interior $\ell^{(o)}$, 
and by $\mu^-$  and $\mu^+$ the end-points of $\ell.$  
The function $\rho_k(\mu)$ maps $\ell^{(o)}$ into the unit circle $U_0$.  
As before, we use the self-adjointness arguments and find $\rho_k'(\mu)\neq0,\;\mu\in\ell^{(o)}.$
Therefore for every $\mu\in\ell^{(o)}$ there exists $t\in[0,2\pi)$ 
such that $\rho(\mu)=e^{it}$ and the local correspondence $t=t(\mu)$ is one-to-one and  
analytic at points of $\ell^{(0) }$. The function ${\rho_k^*({\mu})}$ defined by \eqref{2.14} is 
a single-valued solution of \eqref{1.6} analytic in some neighborhood of $\ell^{(o)}$ and 
hence it coincides with a single-valued branch $\rho_q(\mu)$ of $\rho(\mu)$. Since 
$\rho_k(\mu)\neq\rho_k^*(\mu)$ in a neighborhood of $\mu_0$, we have $q\neq k.$ 

Let now a point $\mu_0$ starts moving monotonically in $\R$ towards 
one of end-points of $\ell$. Then the corresponding 
point $t=t(\mu)$ moves monotonically in $\R$
and simultaneously the points $\rho_k(\mu)$ and $\rho_q(\mu)$ move in the opposite directions 
on the unit circle towards each other. These points may meet for the first time only if 
$\rho_k(\mu)=e^{it}=e^{-it}=\rho_q(\mu)$, i.e., if either $t=0$ or $t=\pi$. For $\mu$ 
corresponding to the meeting point 
we have $\rho_k(\mu)=\rho_q(\mu)=\pm1$ which means that $(\mu,\rho_k(\mu))=(\mu,\rho_q(\mu))$
is a ramification point of $\mathcal{R}(L)$. Therefore
$\mu$ is a point of the set $\mathcal{E}_k(L)$ and since $\ell^{(o)}$ is the maximal open 
interval containing the point $\mu_0$ and no points of $\mathcal{E}_k(L)$, we find that 
$\mu^-\in\mathcal{E}_k(L)$ and $\mu^+\in\mathcal{E}_k(L)$.

The same claim is true if the points $\rho_k(\mu)=e^{it}$ and $\rho_q(\mu)=e^{-it}$ do not 
meet at all. Indeed, in this situation $t\in(0,\pi)$ and because the monotonicity of 
the function
$t(\mu)$ there exist the limits 
$$
t_\pm=\lim\limits_{\mu\to\mu^\pm}t(\mu)\in(0,\pi).
$$
If either $\mu^-$ or $\mu^+$ does not belong to the set $\mathcal{E}_k(L)$, we replace 
$\mu_0$ and $t_0$ in the assumptions of the lemma by either $\mu^-$ and $t_-$ or $\mu^+$ 
and $t_+$, 
respectively, and use the above arguments to find that the interval $\ell$ fails its 
maximal property. The contradiction proves that the end-points of $\ell$ belong to 
 $\mathcal{E}_k(L)$. Differentiating the identity $\rho_k(\mu(t))=e^{it}$
we obtain $\mu'(t)\neq0$ for $t\in \mbox{int}\;S$ where $S$ is the segment with end-points 
$t_-$ and $t_+$, completing the proof of the lemma.

The following statement is a detailed version of
of Theorem \ref{t1}.
\begin{theorem}\label{t3}
For every fixed integer $k=1,...,2n$ there exist a system of non-overlapping compact 
intervals 
$\ell_{j}^{(k)}\subset\sigma_k(L), j=1,...,N(k),$ with end-points in the set 
$\mathcal{E}_k(L)$,
the system of closed intervals $S_{j}^{(k)}, 
j=1,...,N(k),$ each contained in either $[0,\pi]$ or $[\pi,2\pi]$
and the system of continuous monotonic functions $\mu_{jk}(t),$ $t\in S_{j}^{(k)},$ 
$j=1,...,N(k),$ such that

(\;i) The relations 
\begin{equation}\label{2.16}
\Delta(\mu_{jk}(t),\rho_k(\mu_{jk}(t)))=0,\quad\rho_k(\mu_{jk}(t))=e^{it},\quad
t\in{S_{j}^{(k)}}, 
\end{equation}
$$
\mu_{jk}'(t)\neq0,\qquad t\in\mbox{int\;}S_{j}^{(k)},
$$
are valid;

(\;ii) The functions  $\mu_{jk}(t)$ map ${S}_{j}^{(k)}$ one-to-one onto $\ell_{j}^{(k)}$;  

(iii) The representation
$$
\sigma_k(L)=\bigcup\limits_{j=1}^{N(k)}\;{\ell_{j}^{(k)}}
$$
holds;

(\;iv) The spectrum $\sigma(L)$ of operator \eqref{1.3} has the form
$$
\sigma(L)=\bigcup\limits_{k=1}^{2n}\;\left(\bigcup\limits_{j=1}^{N(k)}\;{\ell_{j}^{(k)}}
\right).
$$
\end{theorem}

Proof. Let $\mu_0$ be a non-isolated point in the set $\sigma_k(L)$. If 
$\rho_k(\mu_0)=e^{it_0}$,  $t_0\in[0,\pi]$ and $\mu_0\notin\mathcal{E}_k(L)$ , 
then according to Lemma \ref{l2.2} there exist 
the closed intervals $\ell\in\sigma_k(L)$ and $S\in[0,\pi]$ and the function $\mu(t)$
with properties stated in $(i)$ and $(ii)$. 

If $\rho_k(\mu_0)=e^{it_0}$, $t_0\in[\pi,2\pi]$ and $\mu_0\notin\mathcal{E}_k(L)$ 
we apply Lemma \ref{l2.2} to the function 
$\rho_k^*(\mu_0)$ defined by \eqref{2.14} and to $t_0^*=2\pi-t_0\in[0,\pi]$ in 
the capacity of $\rho_k(\mu_0)$ and $t_0\in[0,\pi]$
and, in addition to the interval $\ell$, find the interval $S^*\in[0,\pi]$ and the 
function $\mu^*(t)$ with properties stated in Lemma \ref{l2.2}. Now the intervals $\ell$ and 
$S=\{t\in[\pi,2\pi]:\;t=2\pi-t^*,\;t^*\in S^*\}$ and the function $\mu(t)=
\mu^*(2\pi-t)$ possess properties  $(i)$ and $(ii)$.

If as before $\mu_0$ is a non-isolated point in the set $\sigma_k(L)$ but 
$\mu_0\in\mathcal{E}_k(L)$, then we fix an interval $V_0=\{\mu\in\R:
|\mu-\mu_0|<\epsilon\}$ without points of the set $\mathcal{E}_k(L)$ distinct from $\mu_0$. 
We choose an arbitrary point $\nu\in V_0\cap\sigma_k(L), \nu\neq\mu_0,$ and again using Lemma \ref{l2.2}
find the interval $\ell$ containing $\nu$, the interval $S$ and the function $\mu(t)$
with properties $(i)$ and $(ii)$. The end-points of $\ell$ must be located in the set 
$\mathcal{E}_k(L)$ and hence $\mu_0$ is one of them.

The remaining opportunity for $\mu_0$ is to be an isolated point in $\sigma_k(L)$. If 
such a point does not belong to $\mathcal{E}_k(L)$, then $\rho_k(\mu)$ is a single-valued 
analytic function in a complex neighborhood of $\mu_0$. As shown in the proof of Lemma
\ref{l2.2}, it maps one-to-one a small interval $\{\mu\in\R:|\mu-\mu_0|<\epsilon\}$ onto a 
small arc of the unit circle $U_0$ centered at $\rho_k(\mu_0)$, and  $\mu_0$ cannot be a 
non-isolated point in $\sigma_k(L)$. Therefore  $\mu_0\in\mathcal{E}_k(L)$ and if 
 $\rho_k(\mu_0)=e^{it_0}$ the conditions $(i)$ and $(ii)$ are satisfied with $\ell=\{\mu_0\}$,
$S=\{t_0\}$ and $\mu(t)\equiv \mu_0$ completing the proof of Theorem \ref{t3}.

The previous analysis shows that the set $\sigma_k(L)$ is formed by a system of compact 
intervals
with end-points in the set $\mathcal{E}_k(L)$. We denote these intervals 
by $\ell_j^{(k)}, j=1,...,N(k),$ 
according to the ordering of their end-points on the real axis and by $S_j^{(k)}$ 
and 
$\mu_{jk}(t)$ the corresponding intervals in $[0,2\pi]$ and functions defined on
them, respectively. 

\begin{remark}\label{r3}
Every point of the set $\mathcal{E}(L)$ belongs to the spectrum $\sigma(L)$ and is 
not isolated in it. Therefore such a point belongs to some non-trivial  interval 
$\ell^{(k)}_j$ and according to Theorem \ref{t3} must be its end-point.
\end{remark}

\section{Operators $L_t$}

For every $t\in[0,2\pi]$ we denote by $L_t$ the selfadjoint operator in the space
$L^2([0,\pi])$ generated by the expression \ref{1.3} and boundary 
conditions 
\begin{equation}
y^{(j)}(\pi)=e^{it}y^{(j)}(0),\qquad 0\leq j\leq 2n-1,\label{3.1}
\end{equation}
The spectrum $\sigma(L_t)$ of $L_t$ in the space ${\cal L}^2[0,\pi]$ coincides with the 
set of all $\mu$'s satisfying Equation \eqref{1.6} with $\rho=e^{it}$ and 
Theorem \ref {t1} states 
$$
\sigma(L)=\bigcup\limits_{t\in[0,2\pi]\;}\sigma(L_t).
$$
In this section we show that the spectrum $\sigma(L_t)$ 
is simple for all $t\in[0,2\pi]$ except maybe 
finitely many points and describe eigen-functions for non-exceptional values of $t.$
\begin{lemma}\label{l3.1}
If $\Delta(\mu_0,e^{it})=0$ and $E(x;\mu,\rho)$ is defined by \eqref{1.12}, then
\begin{equation}\label{3.2}
\|E(.;\mu_0,e^{it})\|^2_{L^2([0,\pi])}=(-1)^{n+1}e^{-it}\Delta'_\mu(\mu_0,e^{it})\;
E(0;\mu_0,e^{-it}).
\end{equation}
\end{lemma}
Proof. 
For every $\mu\in\mathbb{C}$ the definition \eqref{1.12} of $E(x;\mu,\rho)$ 
yields relations
$$
E^{(j)}(\pi;\mu,e^{it})-e^{it}E^{(j)}(0;\mu,e^{it})=0,\qquad j=0,...,2n-2,
$$
$$
E^{(2n-1)}(\pi;\mu,e^{it})-e^{it}E^{(2n-1)}(0;\mu,e^{it})=-\Delta(\mu,e^{it}).
$$
Using the Lagrange formula we find 
$$
(\mu-\mu_0)(E(.;\mu,e^{it}),E(.;\mu_0,e^{it}))=
(-1)^{n+1}e^{-it}\Delta(\mu,e^{it})E(0;\mu_0,e^{-it})
$$
and \eqref{3.2} follows as $\mu\to\mu_0$.

\begin{lemma}\label{l3.2}
The set 
$$
{\mathcal T}(L)=\{t\in[0,2\pi]:
\Delta'_\rho(\mu_0,e^{it})E(0;\mu_0,e^{-it})=0\;\text{for some}\;\mu_0\in\sigma(L_t)\}
$$
is finite and contains the points $0,\pi,2\pi.$
\end{lemma}
Proof. It is evident that the set $\mathcal{T}$ is the union of (maybe intersecting)
sets
$$
{\mathcal T}_1(L)=\{t\in[0,2\pi]:\text{there exists}\; \mu_0\in\sigma(L_t)\; 
\text{such that}\; \Delta'_\rho(\mu_0,e^{-it})=0\}
$$
and
$$
{\mathcal T}_2(L)=\{t\in[0,2\pi]:\text{there exists}\; \mu_0\in\sigma(L_t)\; 
\text{such that}\; E(0;\mu_0,e^{it})=0\}.
$$
If $\mu_0\in\sigma(L_{t_0})$ then there exist the intervals $\ell_{j}^{(k)}$ and 
$S_{j}^{(k)}$ and the function $\mu_{jk}(t),\;t\in S_{j}^{(k)}$ with properties described
in Theorem \ref{t3} such that $\mu_0\in\sigma_k(L),\; \rho_k(\mu_0)=e^{it_0},\;
t_0\in S_{j}^{(k)}$. If $\mu_0\in{\mathcal T}_1(L)$ then $\mu_0\in{\mathcal E_k(L)}$ and 
$\mu_0$ is an end-point of $\ell_{j}^{(k)}$ implying that $t_0$ is an end-point of 
$S_{j}^{(k)}$. Let us show that the set of all $t's$ which are the end-points of intervals 
\begin{equation}\label{2.12}
S_{j}^{(k)},\quad j=1,...,N(k); \quad k=1,...,2n,
\end{equation}
is finite.

For every $R>0$ the number of intervals $\ell_{j}^{(k)}$ either partly or completely  
located in the disc $\{\mu:|\mu|\leq R\}$ is finite and therefore the same is the number 
of intervals $S_{j}^{(k)}$ with the corresponding indices $k$ and $j$. 

It follows from the asymptotic representation \eqref{2.20} that 
there are only two values of $k$ for which $N(k)=\infty$, and they are $n$ and
$n+1$. We choose $R$ large enough for end-points $\alpha_j^{(n)}$  and $\beta_j^{(n)}$ of  
intervals $\ell_{j}^{(n)}$ lying in the domain  $\{\mu:|\mu|\geq R\}$ to be close to the set 
$\{m^{2n}\}_{m\in\N}$. In addition, we can assume that ramification points $(\mu,
\rho_n(\mu))$ of the surface $\mathcal{R}(L)$ with $|\mu|\geq R$ are close to 
$(\alpha_j^{(n)},\rho_n(\alpha_j^{(n)}))$  and $(\beta_j^{(n)},\rho_n(\beta_j^{(n)}))$.

Let  $\ell_{j}^{(n)}$ and  $\ell_{j+1}^{(n)}$ be two adjacent intervals from the set 
$\sigma_n(L)\cap\{\mu:\mu\geq R\}$. If $\gamma_j=(j+1/2)^{2n}$ then it follows from 
\eqref{2.20} that one of the numbers $\rho_n(\gamma_{j+1})$ and $\rho_n(\gamma_j)$ 
is contained in the half-plane $\{\mu:\Im\mu>0\}$ and another in $\{\mu:\Im\mu<0\}$. 
By virtue of Lemma \ref{l2.2} one of intervals $S_{j}^{(n)}$ and $S_{j+1}^{(n)}$ 
belongs to the segment $[0,\pi]$ and another to $[\pi,2\pi]$.

Suppose that the gap $[\beta_j^{(n)},\alpha_{j+1}^{(n)}]$ collapses to a point. If $\mu$
approaches this point then the continuity of  $\rho_n(\mu)$ yields 
$$
\lim\limits_{\mu\to\beta_j^{(n)}-0}\rho_n(\mu)=
\lim\limits_{\mu\to\alpha_{j+1}^{(n)}+0}\rho_n(\mu).
$$
Both numbers here are of the form $e^{it}$, one with $t\in[0,\pi]$ 
and another with $t\in[\pi,2\pi]$,
and we conclude that either $t=0,$ or $t=\pi$ or $t=2\pi$.

If the points $\beta_j^{(n)}$ and $\alpha_{j+1}^{(n)}$ are distinct, then they are separated 
by a non-degenerate open gap not containing points from the set $\sigma_n(L).$ In such a case the 
surface $\mathcal{R}(L)$ is ramified at both points 
$(\alpha_{j+1}^{(n)},\rho_n(\alpha_{j+1}^{(n)}))$  and $(\beta_j^{(n)},\rho_n(\beta_j^{(n)}))$:
otherwise $\rho_n(\mu)$ is a single-valued analytic function in a small neighborhood of each 
such point, representation \eqref{2.13} is valid with $p=1$ and with either 
$\mu_0=\alpha_{j+1}^{(n)}$, or $\mu_0=\beta_{j}^{(n)}$ and hence the gap 
 $(\beta_j^{(n)},\alpha_{j+1}^{(n)})$ contains points of $\sigma_n(L)$ leading us to a 
contradiction. According to the choice of the number $R$ there are only two branches of 
$\rho(\mu)$ in small neighborhoods of points 
$(\alpha_{j+1}^{(n)},\rho_n(\alpha_j^{(n)}))$  and $(\beta_j^{(n)},\rho_n(\beta_j^{(n)}))$, 
and these are $\rho_n(\mu)$ and $\rho^*_n(\mu)=\rho_{n+1}(\mu)$. 
Since 
$$
\lim\limits_{\mu\to\beta_j^{(n)}-0}\;\rho_n(\mu)=
\lim\limits_{\mu\to\beta_{j}^{(n)}-0}\;\rho_n^*(\mu),\quad 
\lim\limits_{\mu\to\alpha_{j+1}^{(n)}+0}\rho_n(\mu)=
\lim\limits_{\mu\to\alpha_{j+1}^{(n)}+0}\rho_n^*(\mu),
$$
$\Im\rho_n(\mu)\Im\rho_n^*(\mu)\leq0$ and all four numbers here are of the form $e^{it,}$ 
they are either $-1$ or $1$ and we again find that either $t=0$ or $t=\pi$ or $t=2\pi.$
Thus the set ${\mathcal T}_1(L)$ is finite.

Let us prove now that the function $E(0;\mu,\rho(\mu))$  does not vanish identically 
on the surface $\mathcal{R}(L)$. 
First we note that this function is one of $(2n)^2$  minors of the matrix $U(\mu)
-\rho(\mu)I$ and hence it is an entry of the matrix
\begin{equation}\label{3.3}
-\lim\limits_{z\to\rho(\mu)}\Delta(\mu,z)(U(\mu)-zI)^{-1}.
\end{equation}
If $\rho(\mu)$ is a simple root of Equation \eqref{1.6} then the previous expression is equal to
\begin{equation}\label{3.4}
\Delta'_\rho(\mu,\rho(\mu))\;res\{(U(\mu)-zI)^{-1};{z=\rho(\mu)}\}.
\end{equation}
Let $\delta>0$ be fixed and let  
$$
S_n(r,\delta)=\{\mu:\mu=\lambda^{2n},\lambda\in T_n(r,\delta)\}\bigcap\{\mu:\Im\mu\geq-1\}
$$ 
where $T_n(r,\delta)$ is defined in Section 2. Then relations \eqref{2.25}
are fulfilled and the estimates 
$$
|e^{\pi i\lambda}-e^{\pi\omega_k\lambda}|\geq d,\qquad \omega_k\neq i,\qquad\lambda\in
T_n(r,\delta)
$$
are valid with a constant $d>0$ not depending on $\lambda$. Therefore 
for every  $\lambda\in T_n(r,\delta)$ and every $z$ satisfying $|z-e^{\pi i\lambda}|=d/4$
we have
$$
|z-e^{\pi\omega_k\lambda}|\geq |e^{\pi i\lambda}-e^{\pi\omega_k\lambda}|-
|z-e^{\pi i\lambda}|\geq3d/4,\quad \omega_k\neq i.
$$
As a result
$$
\|
(E(\pi\lambda)-zI)^{-1}\|\leq C,\quad \lambda\in T_n(r,\delta),
\quad |z-e^{\pi i\lambda}|=d/4,
$$
with a constant $C$ not depending on either $\mu$ or $z$.

To calculate the residue in \eqref{3.4} we again use \eqref{2.25}. Under the same
restrictions  on $\mu$ and $z$ as in the latter relation we have 
$$
(U(\mu)-zI)^{-1}=D(\lambda)
(\Omega+o(1))(E(\pi\lambda)-z(I+o(1))(\Omega^{-1}+o(1))(D(\lambda))^{-1}
$$
with $\mu=\lambda^{2n}, \lambda\in T_n(r,\delta)$ and
$$
\|(E(\pi\lambda)-z(I+o(1))])^{-1}\|=\|(E(\pi\lambda)-zI)^{-1}(I+o(1))\|\leq C.
$$
If the number $r$ is sufficiently large we have the estimates 
$$
|\rho_{n}(\mu)-e^{\pi i\lambda}|\leq d/4,\;|\rho_{n+1}(\mu)-e^{-\pi i\lambda}|\leq d/4,
\; \mu=\lambda^{2n},\;\lambda\in T_n(r,\delta).
$$
Therefore for the same $\mu$'s we obtain
$$
|\rho_k(\mu)-e^{\pi i\lambda}|\geq\left\{
\begin{array}{ccc}
\displaystyle\frac{1}{2}e^{\pi\Re(\omega_k\lambda)}-|e^{\pi i\lambda}|&\geq d,
&1\leq k\leq n-1,\\
 & & \\
\displaystyle|e^{\pi i\lambda}|-\frac{1}{2}e^{\pi\Re(\omega_k\lambda)}&\geq\frac{d}{2},&
n+2\leq k\leq 2n,\\
 & & \\
\displaystyle
d-|\rho_{n+1}(\mu)-e^{-\pi i\lambda}|&\geq\frac{3d}{4},&
k=n+1,
\end{array}
\right.
$$
and hence the eigenvalue $\rho_n(\mu)$ is inside the circle 
$C(\mu,d)=\{z\in\mathbb{C}:|z-e^{\pi i\lambda}|=d/4,\;\mu=\lambda^{2n}\},$
while all other eigenvalues $\rho_k(\mu)$ are outside it. Since 
$$
(E(\pi\lambda)-z(I+o(1)))^{-1}
$$
$$
=(E(\pi\lambda)-zI)^{-1}+
(E(\pi\lambda)-z(I+o(1)))^{-1}o(1)(E(\pi\lambda)-zI)^{-1},
$$
we obtain
$$
res\{(U(\mu)-zI)^{-1};{z=\rho_n(\mu)}\}
$$
$$=\frac{1}{2\pi i}\oint\limits_{C(\mu,d)}\;dz\;
(U(\mu)-zI)^{-1}=D(\lambda)
(\Omega P \Omega^{-1}+o(1))(D(\lambda))^{-1}
$$
where $P=||\delta_{pn}\delta_{jn}||_{j,p=1}^{2n}.$ Thus 
\begin{equation}\label{3.5}
\displaystyle
E(0;\mu,\rho_n(\mu))=\mu^{\frac{2n-1}{2n}}\Delta'_\rho(\mu,\rho_n(\mu))(C+o(1)),
\quad\mu=\lambda^{2n},\quad\lambda\in T_n(r,\delta),
\end{equation}
with a constant $C\neq0$ and  we conclude that $E(0;\mu,\rho(\mu))$  is a 
non-trivial analytic function on the surface $\mathcal{R}(L)$. 
Lemma \ref{l2.0} implies that zeros of this function may accumulate only at the point in infinity, and 
for every $r>0$ there exist finitely many solutions of the equation  $E(0;\mu,\rho(\mu))
=0$ with $|\mu|\leq r$. Since $\rho(\mu)$ is not more than $2n$-valued function, for 
every such solution $\mu$ there exists not more than $2n$ values $t\in(0,2\pi)$
satisfying $\rho(\mu)=e^{it}$.

Furthermore, it follows from \eqref{3.5} that the function $E(0;\mu,\rho_n(\mu))$
may vanish for sufficiently large values of $\Re\mu$ inside exceptional sets 
$D_k=\{\mu=\lambda^{2n}:|\lambda-k|\leq\delta, k\in\N\}$ only. If $E(0;\mu,\rho_n(\mu))=0$
for $\mu\in D_k\cap\sigma_n(L)$ then $E(0;\mu,\rho_{n+1}(\mu))=0$ as well. In the case 
$\rho_n(\mu)=e^{it}\neq\pm1$ we have $\rho_n(\mu)\neq\rho_{n+1}(\mu)$ and 
both $\Delta'_\rho(\mu,\rho_n(\mu))$ and $\Delta'_\rho(\mu,\rho_{n+1}(\mu))$ do not 
vanish. Therefore we can differentiate \eqref{2.16} at the point $\mu$. The resulting identity 
\begin{equation}\label{3.6}
\Delta'_\mu(\mu,e^{it})\mu'_{jk}(t)+ie^{it}\Delta'_\rho(\mu,e^{it})=0
\end{equation}
shows that $\Delta'_\mu(\mu,\rho_n(\mu))$ and $\Delta'_\mu(\mu,\rho_{n+1}(\mu))$ also do 
not vanish. According to \eqref{3.2} the number $\mu$ is at least a double root of both 
$E(0;\mu,\rho_n(\mu))$ and $E(0;\mu,\rho_{n+1}(\mu))$. As a result $\mu$ is a root of 
multiplicity at least $4$ of the function
$E(0;\mu,\rho_n(\mu))E(0;\mu,\rho_{n+1}(\mu))$ which is single-valued and 
analytic inside $D_k$. On the other hand, we have seen in the proof of Theorem \ref{t6}
that for $\mu\in D_k$ and $\rho$ sufficiently close to $\pm1$ the representation 
$\Delta(\mu,\rho)=
(\rho^2-2(\cos(\lambda\pi)+o(1))\rho+1+o(1))\phi(\mu,\rho)$ holds where 
$\phi(\mu,\rho)$ is a polynomial in $\rho$ whose coefficients, as well as all terms $o(1)$,
are single-valued analytic functions of $\mu\in D_k.$ Hence 
$\omega(\mu)\equiv\Delta'_\rho(\mu,\rho_n(\mu))\Delta'_\rho(\mu,\rho_{n+1}(\mu))=
(-4(\sin\pi\lambda)^2+o(1))\psi(\mu)$ is a non-vanishing single-valued analytic function 
of $\mu\in D_k$. We conclude that $\omega(\mu)$ has only two roots in $D_k$ while 
the Rouche Theorem claims, according to \eqref{3.5}, that there are at least 4 such roots. 
Therefore $E(0;\mu,\rho_n(\mu))=0$ is possible for $\rho_n(\mu)=e^{it}=\pm1$ only,
proving that the set ${\mathcal T}_2(L)$ is finite.

\begin{theorem}\label{t4}
If $t\in[0,2\pi]\setminus\mathcal{T}(L)$, then the spectrum $\sigma(L_t)$ of the operator $L_t$ 
is simple and every function $f\in L^2([0,\pi])$ is representable by the 
$L^2([0,\pi])$-convergent orthogonal series
$$
f(x)=\sum\limits_{\Delta(\mu,e^{it})=0}\;w(\mu,e^{it})
E(x;\mu,e^{it})\int\limits_0^\pi dy\;E(y;\mu,e^{-it})f(y)
$$
with the weight function
\begin{equation}\label{3.60}
w(\mu,\rho)=|\Delta'_\mu(\mu,\rho)E(0;\mu,\rho^{-1})|^{-1}.
\end{equation}
\end{theorem}
Proof. For every $\mu\in\sigma(L_{t})$ there exists an integer $k, 1\leq k\leq 2n,$
such that $\mu$ is a point of the set $\sigma_k(L)$ defined by \eqref{1.13}. We apply
Theorem \ref{t3} and find an interval $\ell_{j}^{(k)}$ containing $\mu$, 
an interval $S_{j}^{(k)}\subset[0,2\pi]$ containing $t$ such that $\rho_k(\mu)=e^{it}$
and a function $\mu_{jk}(s)$ satisfying \eqref{2.16}. 
If $t\notin\mathcal{T}(L)$ then $\mu$ does not belong to the set  $\mathcal{E}_k(L)$
defined by \eqref{2.17}. Therefore $\mu$ is an inner point of $\ell_{j}^{(k)}$, 
$t$ is an inner point of $S_{j}^{(k)}$ and \eqref{3.6} 
shows that $\Delta'_\mu(\mu,e^{it})\neq0$. Hence $\mu$ is a simple eigenvalue of 
$L_{t}$, see \cite{NAI}. For $t\notin\mathcal{T}(L)$ we have $E(0;\mu,e^{it})\neq0,$
and it follows from Equation \eqref{3.2} that $E(x;\mu,e^{it})$ is an eigenfunction 
corresponding to $\mu.$ Theorem \ref{t4} states a well-known property of 
selfadjoint operators in Hilbert spaces.

\section{Spectral expansions generated by $L$}

In the present section we prove Theorem \ref{t2}. The main tool in the proof is 
a version of the
Fourier transform proposed by Gel'fand \cite{GEL} for a study of differential operators
with periodic coefficients. 

Given a function $f\in L^2(\mathbb{R})$ we set, following \cite{GEL},
$$
F(x,t)=(\mathcal{G}f)(x,t)=\sum\limits_{r=-\infty}^\infty\;e^{-irt}f(x+\pi r),
\quad t\in[0,2\pi],
$$
and obtain a function $F(x,t)\in L^2(Q)$ where 
$Q=\{(x,t):x\in[0,\pi],t\in[0,2\pi]\}.$ The inverse transform is given by
$$
f(x+\pi r)=\frac{1}{2\pi}\int\limits_0^{2\pi}\;dt\;e^{irt}\;F(x,t),\quad r\in\mathbb{Z},
\quad x\in[0,\pi].
$$

Since
$$
\|F\|^2_{ L^2(Q)}=\frac{1}{2\pi}\iint\limits_{Q}\;dx\;dt\;|F(x,t)|^2
=\int\limits_\mathbb{R}\;dx\;|f(x)|^2=\|f\|^2_{L^2(\mathbb{R})},
$$
the Gel'fand transform $\mathcal{G}$ is an isomorphic linear mapping of $L^2(\mathbb{R})$ onto 
$L^2(Q)$.

For $t\in[0,2\pi]$ and $k=1,...,2n$ we denote by $J_k(t)$ the set of all integers 
$j\in\mathbb{N}$ such that $\ell_{j}^{(k)}\cap\sigma(L_t)\neq\emptyset$. 

Let $0=t_0<t_1<...<t_{M-1}<t_M=2\pi$ be all points of the exceptional set $\mathcal{T}(L)$ 
defined in Lemma \ref{l3.2} and let $t\in(t_m,t_{m+1})$ for some integer $m,\;0\leq m\leq M-1$.
Then to every $j\in J_k(t)$ there corresponds the spectral band $\ell_{j}^{(k)}$, the interval 
$S_{j}^{(k)}\subset[0,2\pi]$ and the function $\mu_{jk}(t), t\in S_{j}^{(k)}$ with properties 
described in Theorem \ref{3.1}. The end-points $m_{jk}$ and of  $M_{jk}$ of 
$S_{j}^{(k)}$ belong to the set $\mathcal{T}(L)$ and 
$$
S_{j}^{(k)}=\bigcup\limits_{p=m_{jk}}^{M_{jk}-1}[t_p,t_{p+1}].
$$
Of course, $m_{jk}\leq m \leq M_{jk}$ and with $s$ moving from $t_m$ to $t_{m+1}$ the
point $\mu_{jk}(s)$ runs over the interval 
$$
\ell_{jm}^{(k)}=\{\mu:\mu=\mu_{jk}(s),s\in(t_m,t_{m+1})\}\subset \ell_{j}^{(k)}
$$
with end-points satisfying either $\mu_{jk}(t_m)<\mu_{jk}(t_{m+1})$ if $\mu_{jk}'(t)>0$ or 
$\mu_{jk}(t_m)>\mu_{jk}(t_{m+1})$ if $\mu_{jk}'(t)<0$ for $t\in S_{j}^{(k)}$. In any case, the set 
$J_k(t)$ does not change with $t$ varying in $(t_m,t_{m+1})$ and 
\begin{equation}\label{4.1}
\ell^{(k)}_j=\bigcup\limits_{p=m_{jk}}^{M_{jk}-1}\ell^{(k)}_{jm}.
\end{equation}

Given an arbitrary function $F\in L^2([0,\pi])$ we set 
$$
\phi(\mu,\rho;F)=\int\limits_0^\pi\;dy\;F(y)\;E(y;\mu,\rho^{-1}).
$$
If here 
$$
\mu=\mu_{jk}(t),\quad\rho=\rho_k(\mu_{jk}(t))=e^{it},\quad 
\quad t\in(t_m,t_{m+1});\quad F=\mathcal{G}f,\quad
f\in L^2(\mathbb{R}),
$$
then 
$$
e^{itr}E(y;\mu_{jk}(t),e^{it})=E(y+\pi r;\mu_{jk}(t),e^{it}),
$$
and hence
\begin{align}
\begin{split}
\phi(\mu_{jk}(t),e^{it};F(.,t))&=\int\limits_0^\pi\;dy\sum\limits
_{r=-\infty}^\infty\;e^{-irt}f(y+\pi r)\;E(y;\mu_{jk}(t),e^{-it})\\
&=\int\limits_0^\pi\;dy\sum\limits_{r=-\infty}^\infty\;f(y+\pi r)
\;E(y+\pi r;\mu_{jk}(t),e^{-it})\\
&=\int\limits_{\mathbb{R}}\;dy\;f(y)\;E(y;\mu_{jk}(t),e^{-it})\\
&=\Phi(\mu_{jk}(t),\rho_{k}(\mu_{jk}(t));f)
\end{split}
\end{align}
where $\Phi(\mu,\rho_k(\mu);f)$ is defined by \eqref{1.7}.

To prove Theorem \ref{t2} assume that  $f$ is an arbitrary function from the 
space $L^2(\mathbb{R})$ with the Gel'fand transform $F=\mathcal{G}f$. According to the 
Fubini Theorem the function $F(.,t)$ belongs to the space $ L^2([0,\pi])$ 
for almost all $t\in[0,2\pi]$, and for every such $t\notin\mathcal{T}(L)$ Theorem \ref{t4}
permits us to represent it by the $ L^2([0,\pi])$-convergent orthogonal series 
\begin{equation}\label{4.2}
F(x,t)=\sum\limits_{k=1}^{2n}\sum\limits_{j\in J_k(t)}
\;
\Psi(\mu_{jk}(t),\rho_k(\mu_{jk}(t));F(.,t))E(x;\mu_{jk}(t),\rho_k(\mu_{jk}(t)))
\end{equation}
with $\Psi(\mu,\rho;f)=w(\mu,\rho)\Phi(\mu,\rho;f)
$ and $w(\mu,\rho)$ defined by \eqref{3.60}.

To prove that the series \eqref{4.2} converges in the space $L^2(Q)$ let us introduce 
the partial sums
$$
F_q(x,t)=\sum\limits_{k=1}^{2n}\;\sum\limits_{j\in J_k(t),\;1\leq j\leq q}
\;
\Psi(\mu_{jk}(t),\rho_k(\mu_{jk}(t));F(.,t))E(x;\mu_{jk}(t),\rho_k(\mu_{jk}(t)))
$$
and set $\Theta(\mu,\rho;f)=w(\mu,\rho)|\phi(\mu,\rho;f)|^2.$ Because the orthogonality 
we have
$$
\|F_q(.,t)\|^2_{L^2([0,\pi])}=\sum\limits_{k=1}^{2n}\;\;
\sum\limits_{j\in J_k(t),\;1\leq j\leq q}
\Theta(\mu_{jk}(t),\rho_k(\mu_{jk}(t);F(.,t)).
$$
Since the set $J_k(t)$ does not depend on $t$ in $(t_m,t_{m+1})$, we can 
integrate the latter sum term by term in every such interval. After substituting $\mu=\mu_{jk}(t)$ into the $j$-th integrated summand we take into account \eqref{4.1} and obtain
$$
\|F_q\|^2_{L^2(Q)}=\frac{1}{2\pi}\sum\limits_{m=1}^{M-1}\int\limits_{t_m}^{t_{m+1}}\;dt
\|F_q(.,t)\|^2_{L^2([0,\pi])}
$$
$$
=\sum\limits_{m=1}^{M-1}\sum\limits_{k=1}^{2n}\sum\limits_{1\leq j\leq q}\;
\int\limits_{\ell^{(k)}_{jm}}\;d\mu\;
p(\mu,\rho_k(\mu))|\Phi(\mu,\rho_k(\mu);f)|^2
$$
$$
=\sum\limits_{k=1}^{2n}\sum\limits_{1\leq j\leq q}\;
\int\limits_{\ell^{(k)}_{j}}\;d\mu\;
p(\mu,\rho_k(\mu))|\Phi(\mu,\rho_k(\mu);f)|^2
$$
According to the Bessel inequality applied to the series \eqref{4.2} 
we have $\|F_q\|^2_{L^2(Q)}\leq\|F\|^2_{L^2(Q)}=\|f\|^2_{\R}$ implying
$$
\sum\limits_{k=1}^{2n}\sum\limits_{j=1}^{N(k)}\;
\int\limits_{\ell^{(k)}_j}\;d\mu\;
p(\mu,\rho_k(\mu))|\Phi(\mu,\rho_k(\mu);f)|^2<\infty.
$$
Since
\begin{equation}\label{4.21}
\|F_{q'}-F_{q''}\|_{L^2(Q)}^2=
\sum\limits_{k=1}^{2n}\sum\limits_{q'<j\leq q''}
\int\limits_{\ell^{(k)}_j}\;d\mu\;
p(\mu,\rho_k(\mu))|\Phi(\mu,\rho_k(\mu);f)|^2,
\end{equation}
the series \eqref{4.2} converges to $F(x,t)$ in $L^2(Q)$-norm and we can 
apply to it the inverse Gel'fand transform term by term. The resulting identity
$$
f(x)=\sum\limits_{k=1}^{2n}\sum\limits_{j=1}^{N(k)}\;
\int\limits_{\ell^{(k)}_j}\;d\mu\;
p(\mu,\rho_k(\mu))\Phi\;(\mu,\rho_k(\mu);f)\;E(x;\mu,\rho_k(\mu))
$$ 
is Equation \eqref{1.8} with $\phi_k(\mu)=\Phi(\mu,\rho_k(\mu);f)$.

Assume now that $\Phi=\{\phi_k(\mu)\}_{k=1}^{2n}\in{\mathcal H}^2_{2n}(L)$. Then
for $k=1,...,2n$ the sequences
$$
\phi_k(\mu)=\{\phi_{jk}(\mu),\; \mu\in\ell^{(k)}_j,\, j=1,2,...\}
$$
satisfy the condition
\begin{equation}\label{4.3}
\sum\limits_{k=1}^{2n}\sum\limits_{j=1}^{N(k)}\;
\int\limits_{\ell^{(k)}_j}\;d\mu\;
p(\mu,\rho_k(\mu))|\phi_{jk}(\mu)|^2<\infty.
\end{equation}
For $t\in(t_m,t_{m+1}), m=0,...,M-1,$ we define
$$
F_q(x,t)=\sum\limits_{k=1}^{2n}\sum\limits_{j\in J_k(t),1\leq j\leq q}
w(\mu_{jk}(t),e^{it})\phi_{jk}(\mu_{jk}(t))\;E(x;\mu_{jk}(t),e^{it})
$$
and obtain
$$
\|F_{q'}(.,t)-F_{q''}(,t)\|_{L^2([0,\pi])}^2=
\sum\limits_{k=1}^{2n}\sum\limits_{j\in J_k(t),\;q'<j\leq q''}
w(\mu_{jk}(t),e^{it})|\phi_{jk}(\mu_{jk}(t))|^2.
$$
Similar to \eqref{4.21} we find
$$
\|F_{q'}-F_{q''}\|_{L^2(Q)}^2=
\sum\limits_{k=1}^{2n}\sum\limits_{q'<j\leq q''}
\int\limits_{\ell^{(k)}_j}\;d\mu\;
p(\mu,\rho_k(\mu))|\phi_{jk}(\mu)|^2.
$$
It follows from \eqref{4.3} 
that $\{F_q\}_{q=1}^\infty$ is a Cauchy sequence 
in the space ${\mathcal H}^2_{2n}(L)$ and hence there exists the function 
$$
F(x,t)=\lim\limits_{q\to\infty}\;F_q(x,t)\in{L^2(Q)}.
$$
Its inverse Gel'fand transform $f=\mathcal{G}^{-1}F$ has the form 
$$
f(x+\pi r)=\frac{1}{2\pi}\int\limits_0^{2\pi}\;dt\;e^{irt}F(x,t)
=\frac{1}{2\pi}\int\limits_0^{2\pi}\;dt\;e^{irt}\lim\limits_{q\to\infty}F_q(x,t)
$$
$$
=\lim\limits_{q\to\infty}\sum\limits_{k=1}^{2n}
\sum\limits_{m=1}^{M-1}\int\limits_{t_m}^{t_{m+1}}dt
\sum\limits_{j\in J_k(t),1\leq j\leq q}
p(\mu_{jk}(t),e^{it})\phi_{jk}(\mu_{jk}(t))\;E(x+\pi r;\mu_{jk}(t),e^{it})
$$
$$
=\sum\limits_{k=1}^{2n}
\sum\limits_{j=1}^{N(k)}
\int\limits_{\ell^{(k)}_j}\;d\mu\;
p(\mu,\rho_k(\mu))\phi_{jk}(\mu)\;E(x+\pi r;\mu,\rho_k(\mu)).
$$
Thus for every $\Phi\in{\mathcal H}^2_{2n}(L)$ there exists a function 
$f\in{L^2(\mathbb{R})}$ such that representation \eqref{1.8} is valid. 
To 
complete the proof of Theorem \ref{t2} let us show that in this representation
$$
\Phi(\mu,\rho_k(\mu);f)=\phi_k(\mu),\qquad \mu\in\sigma_k(L),\qquad k=1,....,2n.
$$
It is sufficient to prove the above relation for an arbitrary step-like function
$$
\Phi_\sigma=\{\phi_k(\mu)\}_{k=1}^{2n},\;\phi_k(\mu)=\left\{
\begin{array}{ccc}
1,&k=k',&\mu\in\sigma\subset\ell^{(k')}_{j'}\\
0& &\text{otherwise}
\end{array}\right.
$$
where $\sigma$ is a closed segment in the interior of a fixed band $\ell^{(k')}_{j'}$
such that the corresponding interval $S^{(k')}_{j'}$ does not contain points of the 
exceptional set 
${\mathcal T}(L)$ and hence the function $w(\mu_{j'k'}(t),e^{it})$ is continuous in it.

First we note that the function $f_\sigma\in{L^2(\mathbb{R})}$ 
defined by \eqref {1.8} with $\Phi=\Phi_\sigma$ has the form
\begin{align}
f_\sigma(x)&=\qquad\int\limits_{\sigma}\;d\mu\;
p(\mu,\rho_{k'}(\mu))\;E(x;\mu,\rho_{k'}(\mu))\\
&=\frac{1}{2\pi}\int\limits_{s^{(k')}_{j'}}\;dt\;
w(\mu_{j'k'}(t),e^{it})\;E(x;\mu_{j'k'}(t),e^{it})
\end{align}
where $s^{(k')}_{j'}\subset S^{(k')}_{j'}$ is the pre-image of $\sigma\subset
\ell^{(k')}_{j'}$ with respect to the function $\mu_{j'k'}(t)$. 

Furthermore, for every integer $k,\; 1\leq k\leq 2n,$ and real 
number $\mu\in\sigma_k(L)\setminus\mathcal{E}_k(L)$ there exist the 
unique integer $j,\; 1\leq j<N(k),$ and real number $\tau\in int\; S^{(k)}_{j}$
such that $\mu=\mu_{jk}(\tau)$. Therefore according to \eqref{1.7}
we obtain
$$
\Phi(\mu,\rho_k(\mu);f_\sigma)=\int\limits_\mathbb{R}\;dy
\;E(y;\mu,\overline{\rho_k(\mu)})f_\sigma(y)
$$
$$
=\lim\limits_{N\to\infty}\sum_{r=-N}^{N}\int_0^\pi\;dy\; 
E(y+\pi r;\mu,\overline{\rho_k(\mu)})f_\sigma(y+\pi r)
$$
$$
=\lim\limits_{N\to\infty}\int\limits_{s^{(k')}_{j'}}dt\;\delta_N(t-\tau)
w(\mu_{j'k'}(t),e^{it})\int_0^\pi dy\;E(y;\mu_{jk}(\tau),e^{-i\tau})
E(y;\mu_{j'k'}(t),e^{it})
$$
where
$$
\delta_N(t)=\frac{1}{2\pi}\sum_{r=-N}^N e^{irt}
$$
is the Dirichlet kernel. For $\mu$ not being an end-point of $\ell^{(k')}_{j'}$ we find
$$
\Phi(\mu,\rho_k(\mu);f_\sigma)
$$
$$=\chi_\sigma(\mu_{jk}(\tau))\;
w(\mu_{j'k'}(\tau),e^{i\tau})\int_0^\pi dy\;E(y;\mu_{jk}(\tau),e^{-i\tau})\;
E(y;\mu_{j'k'}(\tau),e^{i\tau})
$$
where $\chi_\sigma(\lambda)$ is the indicator function of the set $\sigma$. 

If $k\neq k'$ and $\mu_{jk}(\tau)=\mu_{j'k'}(\tau)=\mu$, then $\rho_{k}(\mu)=\rho_{k'}(\mu)$
 contrary to the restriction $\mu\notin{{\mathcal E}_k(L)}$ imposed on $\mu$. 
It means that $\mu_{j'k'}(\tau)\neq\mu_{jk}(\tau)$ for every $k\neq k'$ and the 
eigen-functions $E(y;\mu_{jk}(\tau),e^{i\tau})$ and $E(y;\mu_{j'k'}(\tau),e^{i\tau})$ 
are orthogonal for all natural numbers $j$ and $j'$, yielding 
$$
\Phi(\mu,\rho_k(\mu);f_\sigma)=0,\; k\neq k',\quad k=1,...,2n.
$$
For $k=k'$ we use the orthogonality of  $E(y;\mu_{jk}(\tau),e^{i\tau})$ 
and $E(y;\mu_{j'k}(\tau),e^{i\tau})$ for $j\neq j'$ and \eqref{3.2} and obtain 
$\Phi(\mu,\rho_k(\mu);f_\sigma)=\delta_{kk'}\;\delta_{jj'}\;\chi_\sigma(\mu),\; 
\mu\in\ell^{(k')}_j.$
Therefore $\left\{\Phi(\mu,\rho_k(\mu);f_\sigma)\right\}_{k=1}^{2n}=\Phi_\sigma$
which completes the proof of Theorem \ref{t2}.

\section{Spectral matrix and uniqueness theorem}

Let $\mathcal{M}(\mu)$ be a non-negative Hermitian $2n\times2n$-matrix-valued function
defined on the spectrum $\sigma(L)$. Denote by $L^2_{2n}(\mathcal{M})$ the space 
of complex $2n$-vector functions $F=\{F_q(\mu)\}_{1}^{2n}$ with the scalar product
$$
(F,G)=\int\limits_{\sigma(L)}\;d\mu\;(\mathcal{M}(\mu)F(\mu),G(\mu)).
$$

According to a general definition (cf., \cite{NAI}),
$\mathcal{M}(\mu)$ is called a spectral matrix of $L$ if the relations
\begin{equation}\label{5.1}
F_q(\mu)=\int\limits_\mathbb{R}\;dx\;f(x)u_q(\mu,x),\quad q=1,...,2n,
\end{equation}
and
\begin{equation}\label{5.2}
f(x)=\int\limits_{\sigma(L)}\;d\mu\;(\mathcal{M}(\mu)F(\mu),Y(x,\mu))
\end{equation}
with 
$$
Y(x,\mu)=\text{col}\{u_1(x,\mu),...,u_{2n}(x,\mu)\}
$$
define a one-to-one mapping of the space $L^2(\mathbb{R})$ onto the space 
$L^2_{2n}(\mathcal{M})$ and its inverse, respectively, 
 conjugating operators $L$ and $\mu I$, with 
integrals in \eqref{5.1} and \eqref{5.2} converging in the norms of the corresponding
spaces. As a result, the Parseval identity
\begin{equation}\label{5.3}
\int\limits_\mathbb{R}\;dx\;|f(x)|^2=
\int\limits_{\sigma(L)}\;d\mu\;(\mathcal{M}(\mu)F(\mu),F(\mu))
\end{equation}
holds.

\begin{theorem}\label{t5}
The spectral matrix $\mathcal{M}(\mu)$ of operator $L$ has the form 
\begin{equation}\label{5.4}
\mathcal{M}(\mu)=\sum_{k=1}^{2n}\;\chi_k(\mu)\mathcal{M}(\mu,\rho_k(\mu)),\quad
\mathcal{M}(\mu,\rho)=p(\mu,\rho)||v_q(\mu,\rho)v_{q'}(\mu,\rho^{-1})||_{q,q'=1}^{2n}
\end{equation}
where the function $p(\mu,\rho)$ is defined by \eqref{1.11} and the numbers 
$v_q(\mu,\rho),\;q=1,...,2n$ are uniquely defined by the representation
$$
E(x;\mu,\rho)=\sum\limits_{q=1}^{2n}\;v_q(\mu,\rho)\;u_q(x,\mu).
$$
If two operators of the form \eqref{1.3} have the same spectral matrix 
$\cal M{}$, then their coefficients coincide.
\end{theorem}

Proof. Let $f(x)$ be an arbitrary function from the space $L^2(\mathbb{R})$ and let 
$F_q(\mu)$ be defined by \eqref{5.1}. Then \eqref{1.7} reads
$$
\Phi(\mu,\rho_k(\mu);f)=\sum\limits_{q=1}^{2n}\;v_q(\mu,(\rho_k(\mu))^{-1})\;F_q(\mu)
$$
and Theorem \ref{t2} states that the representation \eqref{5.2} and identity 
\eqref{5.3} hold for functions from the space $L^2(\mathbb{R})$.

On the other hand, let $F=\{F_q(\mu)\}_{q=1}^{2n}$ be an arbitrary element of 
the space $L^2_{2n}(\mathcal{M})$. Then
$$
\int\limits_{\sigma(L)}\;d\mu\;(\mathcal{M}(\mu)F(\mu),F(\mu))<\infty
$$
and if 
$$
\phi_k(\mu)=\sum\limits_{q=1}^{2n}\;v_q(\mu,(\rho_k(\mu))^{-1})\;F_q(\mu)
$$
then the vector $\Phi=\{\phi_k(\mu)\}_{k=1}^{2n}$ belongs to the space 
${\mathcal H}^2_{2n}(L)$. We again use Theorem \ref{t2} and find that the 
function $f(x)$ defined by \eqref{1.8} belongs to the space $L^2(\mathbb{R})$
and the representation itself coincides with \eqref{5.2}. As we have 
just seen, the function $f(x)$ is representable also by \eqref{5.2} with 
$$
\tilde{F}_q(\mu)=\int\limits_\mathbb{R}\;dx\;f(x)u_q(\mu,x),\quad q=1,...,2n.
$$
The Parseval identity \eqref{5.3} implies ${F}=\tilde{F}$ which proves that the mapping 
\eqref{5.1} is from the space $L^2(\mathbb{R})$ onto the space $L^2_{2n}(\mathcal{M})$
and that $\mathcal{M}(\mu)$ is the spectral matrix of $L$.

To prove the second part of Theorem \ref{t5} we assume that the characteristic polynomial
$\Delta(\mu,\rho)$ and the spectral matrix $\mathcal{M}(\mu)$ 
are known. It is easy to see that if $\Delta(\mu,\rho)=0, p(\mu,\rho)\neq\infty,$ and if, 
for a fixed $q'$, the vector 
$$
V(\mu,\rho)=p(\mu,\rho)v_{q'}(\mu,\rho^{-1})\left\{v_q(\mu,\rho)\right\}_{q=1}^{2n}.
$$
does not vanish then it is an eigen-vector of the monodromy matrix $U(\mu)$ 
corresponding to its eigenvalue $\rho$. Now we note that the vector-function 
$V(\mu,\rho(\mu))$ is the $q'$-th line of the matrix $\mathcal{M}(\mu,\rho(\mu))$
and hence it is also known.

Since $v_1(\mu,\rho)=E(0;\mu,\rho)$, we find that $V(\mu,\rho(\mu))$ is a non-trivial
meromorphic vector-function on the surface ${\mathcal R}(L)$ and hence its zeros and poles are 
projected on a discrete set ${\mathcal Z}_V(L)$ of the complex plane 
accumulating at the point in infinity. 

Assume now $(\mu,\rho(\mu))\in{\mathcal R}(L)$ to be such a point that 
$\mu\notin{\mathcal Z}(L)\cup{\mathcal Z}_V(L)$ where ${\mathcal Z}(L)$ 
is a discriminant set of $\Delta(\mu,\rho)$. Then

$\bullet$ The eigenvalues $\rho_1(\mu),...,\rho_{2n}(\mu)$ of the monodromy matrix 
$U(\mu)$ are pairwise distinct;

$\bullet$ The vectors $V(\mu,\rho_1(\mu)),...,V(\mu,\rho_{2n}(\mu))$ are well-defined,
do not vanish and hence form a linear independent system in the space $\R^{2n}$;

$\bullet$ The $2n\times 2n$ matrix $C(\mu,\rho(\mu))$ whose columns are vectors from 
the above system does not degenerate;

$\bullet$ The identity 
$$
U(\mu)=C(\mu,\rho(\mu))\;\text{diag}\{\rho_1(\mu),...,\rho_{2n}(\mu)\;\}
C(\mu,\rho(\mu))^{-1},
\;\mu\notin{\mathcal Z}(L)\cup{\mathcal Z}_V(L),
$$
holds.

The factors in the above product are neither single-valued no analytic in $\mathbb{C}$: 
every rotation of a point $\mu$ around the projection of a ramification point results 
in a permutation of columns and at every point from the set ${\mathcal Z}_V(L)$ the matrix 
$C(\mu,\rho(\mu))^{-1}$ may have a pole. Nevertheless, the product is single-valued in 
$\mathbb{C}\setminus({\mathcal Z}(L)\cup{\mathcal Z}_V(L))$ and coincides with 
the entire matrix-function $U(\mu)$ 
outside the set of all its singular points. Hence all these points are removable 
singularities of the product which means that the spectral matrix  $\mathcal{M}(\mu)$ 
permits us to reconstruct completely the monodromy matrix. According to a 
uniqueness theorem proved by Leibenzon \cite{LEI} the latter uniquely determines all 
coefficients of operator \eqref{1.3} which completes the proof of Theorem \ref{t5}.

It follows from Theorem \ref{t5} that 
the spectral matrix $\mathcal{M}(\mu)$ determines the coefficients of operator \eqref{1.3}
uniquely. For Hill's operators \eqref{1.1} on the entire real axis  
 this statement is well-known and is a version of 
the uniqueness theorem proved by Marchenko \cite{MAR} for Sturm-Liouville operators on a 
semi-axis and extended by Rofe-Beketov \cite{RFB} to such operators on the entire real axis.

\ Department of Mathematics 

\ Ben-Gurion University of the Negev, 

\ Beer--Sheva 84105, Israel
\bigskip

 AMS classification \# 34B05, 34L05

\end{document}